\documentclass[12pt]{amsart}

\usepackage[dvips]{graphicx}
\usepackage{amsfonts,amssymb,amscd,bm}
\newcommand{\R}{{\mathbb R}}

\newcommand{\Z}{{\mathbb Z}}

\newcommand{\T}{{\mathbb T}}

\newcommand{\ph}{{\varphi}}

\newcommand{\be}{{\beta}}
\newcommand{\la}{{\lambda}}

\newcommand{\ga}{{\gamma}}

\newcommand{\De}{{\Delta}}
\newcommand{\si}{{\sigma}}

\newcommand{\Om}{{\Omega}}

\newcommand{\bu}{\bullet}

\newcommand{\pr}{{\mathrm{pr}}}

\newcommand{\lk}{{\mathrm{lk}}}
\newcommand{\tg}{{\mathrm{TG}}}
\newcommand{\fqs}{{\mathrm{fqs}}}
\newcommand{\fes}{{\mathrm{fes}}}

\newcommand{\crit}{ \begin{picture}(10,10)(0,0)
  \put(5,0){\oval(8,8)[t]}
  \put(5,0){\line(0,1){8}}
 \end{picture} }

\newcommand{\cubicsn}{ \begin{picture}(8,10)(0,2)
  \put(3,0){\line(0,1){10}}
  \put(-1,5){\oval(8,8)[br]}
  \put(7,5){\oval(8,8)[tl]}
 \end{picture} }

\newcommand{\cuspint}{ \begin{picture}(12,10)(0,2)
  \put(0,1){\line(2,1){10}}
  \put(1,4){\oval(8,10)[tr]}
  \put(9,4){\oval(8,10)[tl]}
 \end{picture} }
\newcommand{\degcusp}{ \begin{picture}(12,10)(0,2)
  \qbezier(0,10)(10,10)(10,0)
  \qbezier(0,6)(8,6)(10,0)
 \end{picture} }

\newcommand{\horcusp}{ \begin{picture}(8,10)(0,2)
  \put(0,9){\oval(14,8)[br]}
  \put(0,1){\oval(14,8)[tr]}
 \end{picture} }
\newcommand{\hortang}{ \begin{picture}(10,10)(0,2)
  \put(5,1){\oval(8,8)[t]}
  \put(5,9){\oval(8,8)[b]}
 \end{picture} }
\newcommand{\hortrip}{ \begin{picture}(12,10)(0,2)
  \put(5,1){\oval(8,8)[t]}
  \qbezier(3,1)(5,5)(7,9)
  \qbezier(7,1)(5,5)(3,9)
 \end{picture} }
\newcommand{\inctrip}{ \begin{picture}(16,12)(0,2)
  \qbezier(-1,0)(7,14)(15,0)
  \qbezier(5,1)(7,7)(9,13)
  \qbezier(2,1)(7,6)(12,12)
 \end{picture} }

\newcommand{\maxtang}{ \begin{picture}(15,10)(0,0)
  \qbezier(0,0)(6,12)(12,0)
  \qbezier(3,0)(6,12)(9,0)
 \end{picture} }
\newcommand{\quadrup}{ \begin{picture}(12,10)(0,2)
  \put(3,0){\line(2,5){4}}
  \put(0,3){\line(5,2){10}}
  \put(7,0){\line(-2,5){4}}
  \put(0,7){\line(5,-2){10}}
 \end{picture} }

\newcommand{\tangint}{ \begin{picture}(12,10)(0,2)
  \put(0,3){\line(2,1){10}}
  \put(1,5){\oval(8,8)[r]}
  \put(9,5){\oval(8,8)[l]}
 \end{picture} }

\newtheorem{theorem}{Theorem}[section]

\newtheorem{proposition}[theorem]{Proposition}
\newtheorem{lemma}[theorem]{Lemma}

\theoremstyle{definition}
\newtheorem{definition}[theorem]{Definition}
\newtheorem{example}[theorem]{Example}

\title{Fiber quadrisecants in knot isotopies}

\author[fiedler]{T.~Fiedler}
\address{ Laboratoire Emile Picard, Universit\'e Paul Sabatier,
118 route Narbonne, 31062 Toulouse, France}
\email{ fiedler@picard.ups-tlse.fr }

\author[kurlin]{V.~Kurlin}
\address{ Department of Mathematical Sciences,
University of Liverpool,
Liverpool L69 7ZL, United Kingdom}
\email{ kurlin@liv.ac.uk, vak26@yandex.ru }

\date{ January 30, 2007} 

\begin{document}


\begin{abstract}
Fix a straight line $L$ in Euclidean 3-space and consider
 the fibration of the complement of $L$ by half-planes.
A generic knot $K$ in the complement of $L$ has neither
 fiber quadrisecants nor fiber extreme secants such that 
 $K$ touches the corresponding half-plane at 2 points.
Both types of secants occur in generic isotopies of knots.
We give lower bounds for the number of these fiber secants
 in all isotopies connecting given isotopic knots.
The bounds are expressed in terms of invariants calculable
 in linear time with respect to the number of crossings.
\end{abstract}

\maketitle



\section{Introduction}

In this paper we give another application of 
 the main result of \cite{FK}, namely 
 the \emph{higher order} Reidemeister theorem
 for one-parameter families of knots.
Fix a straight line $L$ in $\R^3$, the \emph{axis}.
For simplicity assume that $L$ is horizontal.
Consider the fibration $\ph:\R^3-L\to S_{\ph}^1$ by half-planes
 attached to the axis $L$.
The fibration $\ph$ can be visualized as an open book 
 whose half-planes are fibers of $\ph$.
We will study some distances between isotopic knots in the complement $\R^3-L$.
\smallskip

A knot is the image of a $C^{\infty}$-smooth embedding $S^1\to\R^3-L$.
An isotopy of knots is a smooth family $\{K_t\}$, $t\in[0,1]$, of smooth knots.
The theory of knots in $\R^3-L$ covers the classical knot theory 
 in $\R^3$ and closed braids.
An $n$-braid $\be$ is a family of $n$ disjoint strands in a vertical cylinder 
 such that the strands have fixed enpoints on the horizontal bases of the cylinder 
 and they are monotonic in the vertical direction.
After identifying the bases of the cylinder in Fig.~1,
 any braid $\be$ converts into the \emph{closed} braid $\hat\be$,
 a link in a solid torus going around the axis $L$.
The boundary circle of the lower base of the cylinder plays 
 the role of $L\cup\infty$.

\begin{figure}[!h]
\includegraphics[scale=1.0]{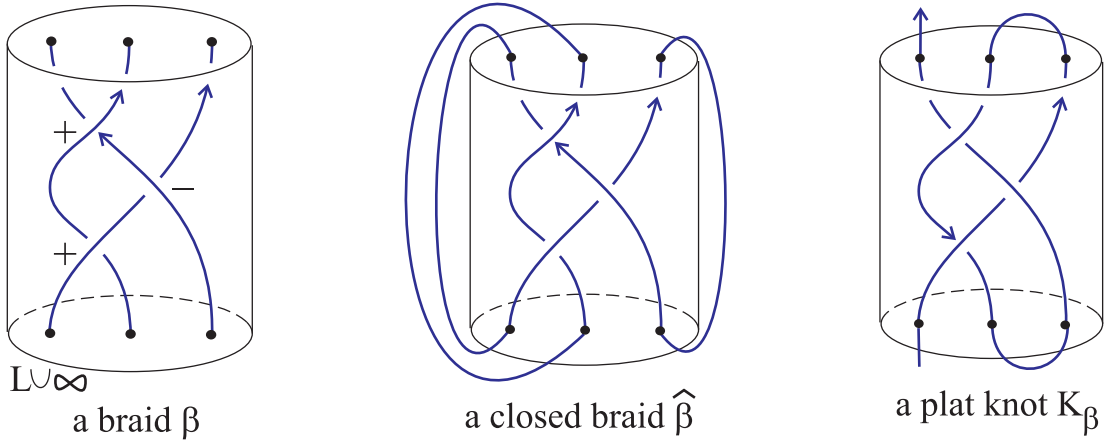}
\caption{Examples of a braid, a closed braid, a plat knot.}
\end{figure}

A \emph{secant}, a \emph{trisecant} and a \emph{quadrisecant} of $K\subset\R^3-L$ is 
 a straight line meeting $K$ transversally in 
 2, 3 and 4 points, respectively.
A secant meeting $K$ in points $p,q$ is \emph{extreme} if 
 the secant and the tangents of $K$ at $p,q$ lie in a common plane.
Namely, $K$ has tangencies of \emph{order} 1 at $p,q$ 
 with a plane passing through the secant, i.e. the plane and $K$
 are given by $\{z=0\}$ and $\{y=0,\, z=x^2\}$ in some local coordinates near $p,q$.
A generic knot has finitely many extreme secants and quadrisecants.
If we are interested only in fiber secants respecting $\ph$ then 
 these geometric features define 
 codimension~1 singularities in the space of all smooth knots $K\subset\R^3-L$.

\begin{definition}
A \emph{fiber} secant, a \emph{fiber} trisecant, 
 a \emph{fiber} quadrisecant of a knot $K\subset\R^3-L$
 is a straight line meeting $K$ transversally in 2, 3, 4 points, respectively, 
 that lie in a fiber of the fibration $\ph:\R^3-L\to S_{\ph}^1$.
A fiber secant meeting $K$ in points $p,q$ is called \emph{extreme}
 if $K$ has tangencies of order~1 at $p,q$ with the fiber.
\end{definition}

We use fiber secants to measure a distance 
 between different embeddings of a knot.
A similar distance with respect to
 Reidemeister moves of type III was studied in \cite{CESS}, see Fig.~4.
Reidemeister moves can be performed on a knot $K$
 in a small neighbourhood of a disk.
Reidemeister moves III correspond to 
 triple points in the horizontal disk of a projection, i.e. 
 to vertical trisecants meeting $K$ in 3 points.
\smallskip

So the authors of \cite{CESS} found the minimal number of vertical
 trisecants in isotopies between different representations of a knot.
We consider more general features of a knot, namely quadrisecants
 in the half-planes of the fibration $\ph$ and estimate
 their minimal number in knot isotopies.
Arbitrary quadrisecants provide lower bounds for
 the ropelength of knots \cite{DDS}.
To define our lower bounds we associate to each knot $K\subset\R^3-L$ 
 an oriented graph $\tg(K)$, the \emph{trace graph} in a thickened torus.
\smallskip

Choose cylindrical coordinates $\rho,\ph,\la$ in $\R^3$, where
 $\la$ is the coordinate in the oriented axis $L$,
 $\rho$ and $\ph$ are polar coordinates in a plane orthogonal to $L$.
For an ordered pair of points $(p,q)\subset\{\ph=\mbox{const}\}$, 
 let $\tau(p,q)$ be the angle between $L$ and the oriented line $S(p,q)$ 
 passing first through $p$ and after through $q$.
Denote by $\rho(p,q)$ the distance between $S(p,q)$
 and the origin $0\in L$.
Introduce the oriented thickened torus 
 $\T=S_{\tau}^1\times S_{\ph}^1\times\R^+_{\rho}$
 parametrized by $\tau,\ph\in[0,2\pi)$ and $\rho\in\R^+$.

\begin{definition}
Take a knot $K\subset\R^3-L$ in general position such that
 $K$ intersects each fiber of $\ph$ in finitely many points.
Map an ordered pair $(p,q)\subset K\cap\{\ph=\mbox{const}\}$
 to $(\tau(p,q),\ph,\rho(p,q))\in\T$.
So each oriented fiber secant of $K$ maps to 
 a point in the thickened torus $\T$.
The image of this map is the \emph{trace graph} $\tg(K)\subset\T$.
\end{definition}

\begin{figure}[!h]
\includegraphics[scale=1.0]{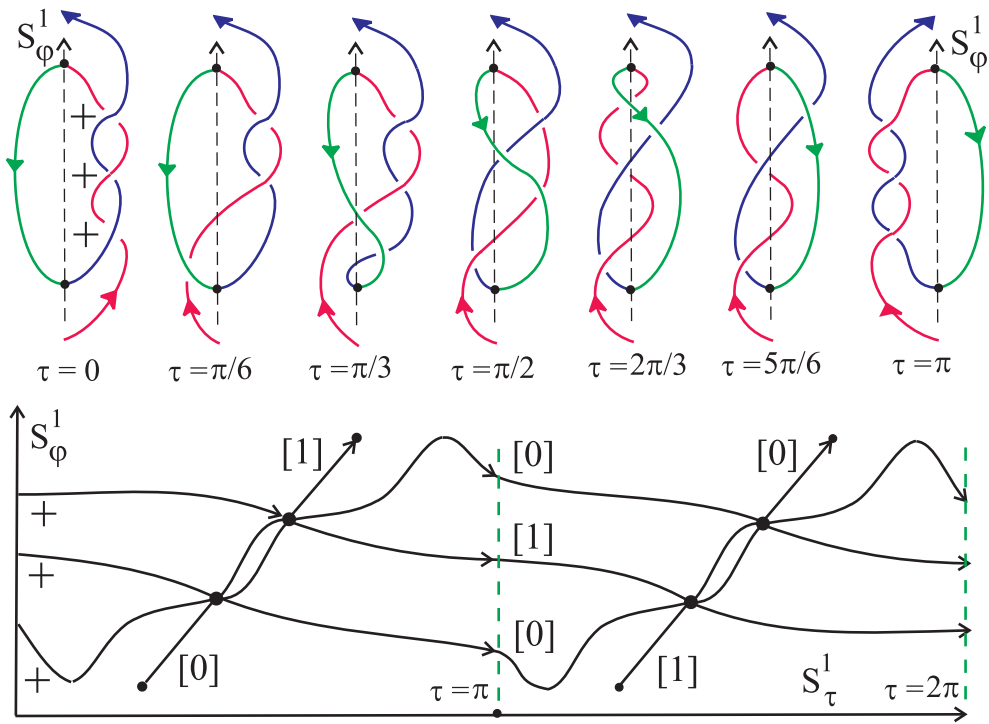}
\caption{The trace graph $\tg(K)$ of the long trefoil $K$.}
\end{figure}

Figure~2 shows the trace graph $\tg(K)$ of a long trefoil $K$ 
 going once around a very long circle $L\cup\infty$.
The fibers there are horizontal planes.
The knots in Fig.~2 are obtained from $K$ by the rotation
 around a vertical line.
A crossing in the projection of a rotated knot corresponds
 to a fiber secant of $K$, i.e. to a point of $\tg(K)$.
The embedding $\tg(K)\subset\T$ is symmetric under 
 the shift $\tau\mapsto\tau+\pi$.
\smallskip

For a generic knot $K$ of Definition~2.1, $\tg(K)$ can have
 only \emph{hanging} vertices and \emph{triple} vertices
 associated to fiber tangents and fiber trisecants of $K$, respectively.
A double crossing of $\tg(K)$ under
 $\pr_{\tau\ph}:\tg(K)\to S_{\tau}^1\times S_{\ph}^1$
 corresponds to a pair of parallel secants meeting $K$ in points
 that lie in a fiber $\{\ph=\mbox{const}\}$.
\smallskip

Let $m$ be the linking number of a knot $K$ with the axis $L$. 
It turns out that the trace graph $\tg(K)$ splits 
 into a union of oriented \emph{traces}
 (arcs or circles) marked by canonically defined homological markings in $\Z_{|m|}$, 
 where $\Z_0=\Z$ and $\Z_1=\{0\}$, see Definition~2.2. 
For example, the closure of $\si_3\si_2\si_1\in B_4$ has the trace graph
 in Fig.~3, which is a disjoint union of 3 trace circles marked by 
 $[1],[2],[3]\in\Z_4$.

\begin{figure}[!h]
\includegraphics[scale=1.0]{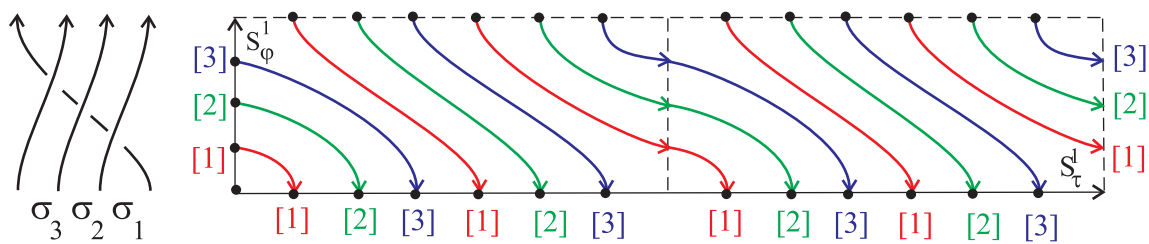}
\caption{The trace graph of $\widehat{\si_3\si_2\si_1}$ splits into 3 trace circles.}
\end{figure}

Introduce the \emph{sign} of a crossing in the projection 
 $\pr_{\tau\ph}(\tg(K))$ as usual, see Fig.~1.
We shall define 3 functions on $\tg(K)$, which will be invariant
 under regular isotopy of $\tg(K)$, not allowing 
 Reidemeiser moves of type I, see Lemma~3.1.

\begin{definition}
Take a knot $K\subset\R^3-L$ such that $\lk(K,L)=m\neq\pm 1$ and
 the projection $\pr_{\tau\ph}(\tg(K))$ has finitely many crossings.
For distinct $[a],[b]\in\Z_{|m|}-\{0\}$, the \emph{unordered} writhe $W_{a,b}^u(K)$ 
 is the sum of signs over all crossings of 
 the trace marked by $[a]$ with the trace marked by $[b]$.
The \emph{ordered} writhe $W_{a,b}^o(K)$ is the sum of signs over all crossings, 
 where the trace $[a]$ overcrosses the trace $[b]$.
The \emph{coordinated} writhe $W_{a,a}^c(K)$ is the sum of signs 
 over all self-crossings of the trace $[a]$.
\end{definition}

We do not consider knots $K\subset\R^3-L$ with $\lk(K,L)=\pm 1$, 
 because in this case $\tg(K)$ splits into trace arcs marked by 
 $[0]$ and $[1]$ only, see Definition~2.2.
The trace graph $\tg(K)$ is constructed
 from a plane projection of $K$, see Lemma~2.4.
The writhes of Definition~1.3 depend on a geometric embedding 
 $K\subset\R^3-L$, but can be computed in linear time with respect
 to the number of crossings of $K$ and change under 
 knot isotopies in a controllable way, see Lemma~3.2.

\begin{theorem}
For isotopic generic knots $K_0,K_1\subset\R^3-L$,
 denote by $\fqs(K_0,K_1)$ and $\fes(K_0,K_1)$ the minimum number of 
 fiber quadrisecants and fiber extreme secants, respectively,
 occuring during all isotopies $\{K_t\}$, $t\in[0,1]$.
\smallskip

\noindent
For isotopic knots $K_0,K_1$, we have
 $\fqs(K_0,K_1)\geq \dfrac{1}{12}\sum\limits_{0\neq a\neq b\neq 0} 
 |W_{a,b}^u(K_0)-W_{a,b}^u(K_1)|$\\
and
 $\fqs(K_0,K_1)+\dfrac{1}{6}\fes(K_0,K_1)\geq \dfrac{1}{12}\sum\limits_{a\neq 0} 
 |W_{a,a}^c(K_0)-W_{a,a}^c(K_1)|$.
Given isotopic closed braids $\hat\be_0,\hat\be_1$, we get
 $\fqs(\hat\be_0,\hat\be_1)\geq \dfrac{1}{12}\sum\limits_{0\neq a\neq b\neq 0} 
 |W_{a,b}^o(\hat\be_0)-W_{a,b}^o(\hat\be_1)|$.
\end{theorem}

The third lower bound is not less than the first one,
 but works for closed braids only.
The second bound gives another estimate
 for the number of fiber quadrisecants for closed braids
 since fiber extreme secants do not occur in braid isotopies.


\section{The trace graph of a knot}

We shall define generic knots $K\subset\R^3-L$ and 
 geometric features of knots, considered as codimension~1 singularities 
 in the space of all knots in $\R^3-L$.
Each singularity is illustrated by a small portion of 
 the projection of $K$ along the corresponding secant.
For example, a tangent of $K$ maps to a cusp
 in the plane projection of $K$ along the tangent, 
 while a quadrisecant projects to a quadruple point. 

\begin{definition}
A knot $K\subset\R^3-L$ is \emph{generic} 
 if $K$ has no following features: 
\smallskip

\noindent
$\quadrup$ : a fiber quadrisecant intersecting $K$
 transversally in 4 points;
\smallskip

\noindent
$\tangint$ : a fiber trisecant meeting $K$ in 3 points such that
 the secant lies in the plane spanned by the tangents of $K$ 
 at 2 of these points;
\smallskip

\noindent
$\cuspint$ : a fiber secant meeting $K$ in 2 points and
 having a tangency of order~1 with $K$ at one of these points;
\smallskip

\noindent
$\cubicsn$ : a fiber secant meeting $K$ in points $p,q$ such that
 $K$ has a tangency of \emph{order}~2 at $p$ with the plane 
 spanned by the secant and the tangent of $K$ at $q$, 
 i.e. the plane and $K$ are given by $\{z=0\}$ and $\{y=0,\, z=x^3\}$
 in local coordinates near $p$;
\smallskip

\noindent
$\degcusp$ : a fiber tangent having a tangency of \emph{order}~2 with $K$,
 i.e. the tangent and $K$ are given by $\{y=z=0\}$ and $\{y=0,\, x^2=z^5\}$
 in local coordinates;
\smallskip

\noindent
$\hortrip,\inctrip$ : a fiber trisecant meeting $K$ in 3 points such that
 $K$ has a tangency of order~1 with the fiber at one of these points;
\smallskip

\noindent
$\horcusp$ : a fiber tangent meeting $K$ in a point, 
 where $K$ has a tangency of \emph{order}~2 with the fiber, i.e. 
 the fiber and $K$ are given locally by $\{z=0\}$ and $\{y=0,\, z=x^3\}$;
\smallskip

\noindent
$\hortang,\maxtang$ : a fiber secant meeting $K$ in 2 points, 
 where $K$ has tangencies of order~1 with the fiber.
\end{definition}

The singularities of Definition~2.1 can be visualized by rotating 
 a knot around a vertical axis.
The last singularity represents 
 two local extrema with same vertical coordinate:
 they collide under the projection after rotating by a suitable angle.
\smallskip

A \emph{trace} in the trace graph $\tg(K)$ of a knot $K$
 is either a subarc ending at hanging vertices or 
 a subcircle of $\tg(K)$.
A trace passes through triple vertices without changing its direction.
The trace graph in Fig.~2 consists of 2 trace arcs.
\smallskip

By Definition~1.2 any point in $\tg(K)$
 corresponds to a fiber secant of $K$ and also to
 an intersection in the projection of $K$ along the secant.
Recall that hanging vertices and triple vertices of $\tg(K)$
 correspond to cusps and triple intersections. 
Mark also \emph{tangent} vertices of degree~2 in $\tg(K)$ such that
 the corresponding secant projects to a tangent point of order~1,
 see 2 tangent vertices in Fig.~6iv.
\smallskip

All points of $\tg(K)$ apart from the vertices of $\tg(K)$ 
 correspond to double crossings with well-defined signs.
Associate to such a point 
 the \emph{sign} of the corresponding crossing in the projection.
While we travel along a trace of $\tg(K)$ the sign
 does not change at triple vertices, but switches at tangent vertices.
\smallskip

Let us look at the function $\tau$ on fiber secants passing
 through 2 points $p,q\in K$.
Namely, $\tau(p,q)$ is the angle between $L$ and the fiber secant through $p,q$.
The function $\tau(p,q)$ has a local extremum if and only if 
 the corresponding secant of $K$ projects to a tangent point,
 i.e. $\tau$ changes its monotonic type at tangent vertices of $\tg(K)$.

\begin{definition}
Take a generic knot $K\subset\R^3-L$ with $\lk(K,L)=m$.
Split $\tg(K)$ by tangent vertices into arcs
 with associated signs coming from plane projections.
Orient each arc so that if the angle $\tau$ is increasing 
 (respectively, decreasing) along the arc then 
 the associated sign of the arc is +1 (respectively, $-1$), see Fig.~2.
\smallskip

Any point of $\tg(K)$ apart from the vertices of $\tg(K)$ 
 is associated to a crossing $(p,q)$ in the projection of $K$
 along the secant through $p,q\in K$.
Smoothing the projection at $(p,q)$ produces a 2-component link.
The linking number of $L$ with the component, 
 where the undercrossing goes to the overcrossing,
 is called the \emph{homological marking} $[a]\in\Z_{|m|}$
 of $(p,q)$ and of the point of $\tg(K)$, see Fig.~3.
\end{definition}

The trace graph in Fig.~2 splits into 2 trace arcs marked by $[0]$ and $[1]$.
Under the shift $\tau\mapsto\tau+\pi$, the homological marking $[a]$
 converts into $[|m|-a]\in\Z_{|m|}$, see Fig.~3.
Recall that a hanging vertex of $\tg(K)$ corresponds to a fiber tangent of $K$,
 i.e. to an ordinary cusp in the plane projection along this tangent.

\begin{figure}[!h]
\includegraphics[scale=1.0]{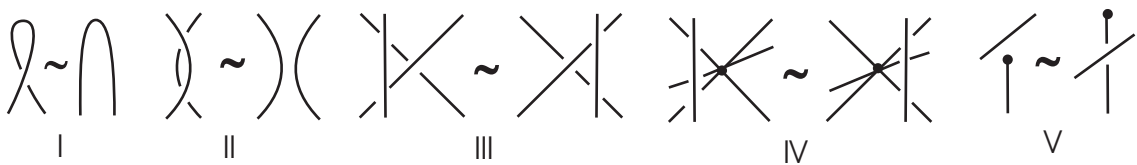}
\caption{Reidemeister moves for trace graphs.}
\end{figure}

\begin{lemma}
The trace graph $\tg(K)$ of a generic knot $K$
 splits into traces with well-defined homological markings.
The orientation of edges, introduced in Definition~2.2,
 provides orientations of all traces of $\tg(K)$.
\end{lemma}
\begin{proof}
Consider the sign of a crossing, monotonic type of the function $\tau$
 and homological marking as functions of a point in the trace graph $\tg(K)$.
All these functions remain constant while the projection of $K$
 along the corresponding secant keeps its combinatorial type.
By the classical Reidemeister theorem, a knot projection can change
 under Reidemeister moves of types I, II, III, see Fig.~4.
\smallskip

Under Reidemeister move I a fiber secant of $K$
 appears or disappears, i.e. the corresponding point in the trace graph
 comes to a hanging vertex of $\tg(K)$.
Under Reidemeister move II, two crossings with 
 opposite signs and same marking appear or disappear.
At this moment the function $\tau$ reverses its monotonic type.
So the orientations of adjacent arcs of $\tg(K)$ agree at tangent vertices.
Under Reidemeister move III nothing changes, i.e.
 all arcs of a trace have same marking.
\end{proof}

The right picture in Fig.~1 shows a \emph{plat} diagram of 
 a knot $K_{\be}$ associated to a braid $\be$.
Any knot can be isotoped to a curve with a plat diagram.

\begin{lemma}
Let $K_{\be}$ be a knot with a plat diagram associated to 
 a $(2n+1)$-braid $\be$ of length $l$.
The trace graph $\tg(K_{\be})$ can be constructed combinatorially
 from the diagram of $K_{\be}$.
The writhes of Definition~1.3 can be computed with complexity $Cln^2$.
\end{lemma}
\begin{proof}
Describe the trace graphs of elementary braids
 containing one crossing only.
Figure 5 shows the explicit example for
 the crossing $\si_1$ of first two strands in the 4-braid.
Firstly we draw all strands in a vertical cylinder.
Secondly we approximate with the first derivative
 the strands forming a crossing by smooth arcs.

\begin{figure}[!h]
\includegraphics[scale=1.0]{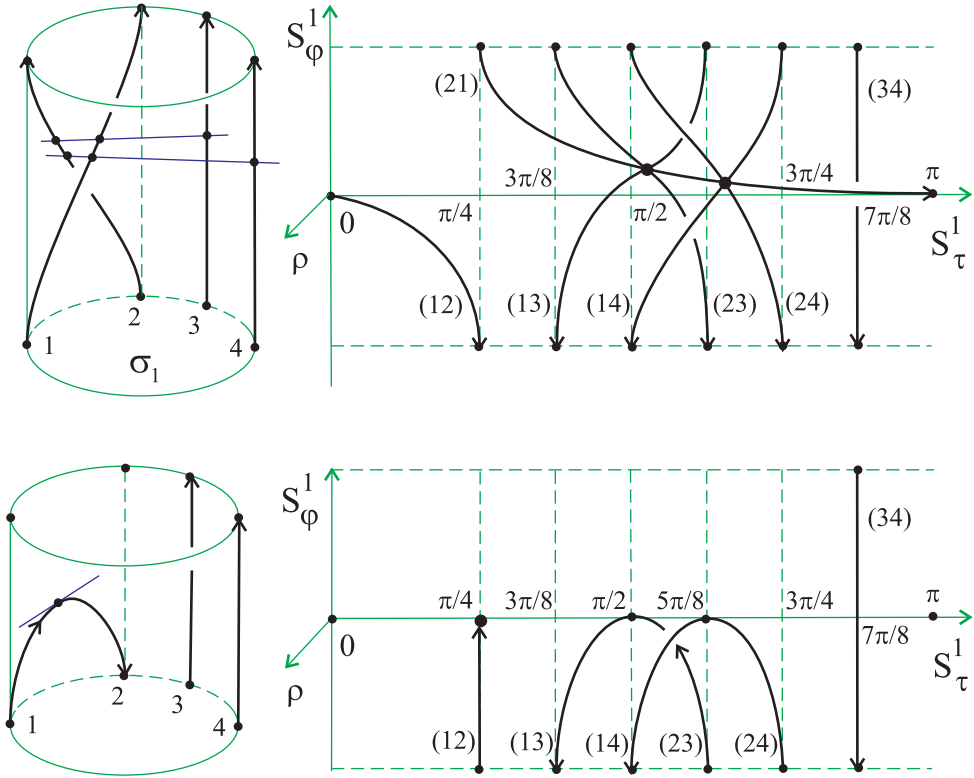}
\caption{Half trace graphs of $\si_1\in B_4$ and a local maximum.}
\end{figure}

The monotonic strands on 
 the left pictures in Fig.~5 are denoted by 1,2,3,4.
The trace graphs on the right pictures have arcs labelled
 by ordered pairs $(ij)$, $i,j\in\{1,2,3,4\}$.
The arc $(ij)$ represents crossings, where the $i$th strand
 overcrosses the $j$th one.
For instance, at the moment $\tau=0$ the braid $\si_1$ has
 exactly one crossing $(12)$, which becomes a crossing $(21)$
 after rotating the braid by $\tau=\pi/4$.
Two triple vertices on the upper right picture correspond to      
 two horizontal trisecants on the upper left picture in Fig.~5.
Similarly we construct the trace graph of a local extremum.
The only hanging vertex corresponds to a horizontal tangent.

In general we split $K_{\be}$ by fibers of $\ph:\R^3-L\to S_{\ph}^1$
 into several sectors each of that contains exactly
 one crossing or one extremum.
To each sector we associate the corresponding
 elementary block and glue them together.
The resulting trace graph contains $2l(2n-1)$ triple vertices
 and $2n$ hanging vertices.
Any two arcs in an elementary block have at most one crossing,
 i.e. not more than $n^2$ crossings in total.
Hence each writhe of Definition~1.3 can be computed
 with complexity $Cln^2$.
\end{proof}

\begin{definition}
Denote by $\Om$ the discriminant of knots $K$ failing to be generic
 due to one of the singularities of Definition~2.1.
An isotopy of knots $\{K_t\}$, $t\in[0,1]$, is \emph{generic} if
 the path $\{K_t\}$ intersects $\Om$ transversally.
A \emph{regular} isotopy of trace graphs is generated by 
 the Reidemeister moves of types II, III, IV, V in Fig.~4.
\end{definition}

Any orientations and symmetric images of the moves in Fig.~4 are allowed.
Proposition~2.6 is a particular case of 
 a more general higher order Reidemeister theorem \cite[Theorem~1.8]{FK}.
A knot can be reconstructed from its trace graphs equipped with
 labels, ordered pairs of integers, see details in \cite[section~5]{FK}.

\begin{proposition}
If knots $K_0,K_1\subset\R^3-L$ are isotopic then
 $\tg(K_0),\tg(K_1)$ are related by 
 regular isotopy and a finite sequence of the moves in Fig.~6.
\end{proposition}
\begin{proof}
The singularities of Definition~2.4 are all codimension~1 singularities
 associated to fiber secants and fiber tangents of knots, 
 see \cite[section~3]{FK}.
Any isotopy of knots can be approximated by a generic isotopy of Definition~2.5.
Each move of Fig.~6 corresponds to one of the singularities.
For instance,  when a path in the space of knots passes through 
 a knot with a fiber quadrisecant, the tetrahedral move 6i changes 
 the trace graph by collapsing and blowing up a tetrahedron.
A formal correspondence between the singularities and moves
 was shown in \cite[Claim~4.5]{FK}.
Empty vertices of degree~2 in Fig.~6 denote points corresponding
 to the singularity $\crit$, where a knot touches a fiber and has
 a fiber secant through this tangent point.
\smallskip

The moves of Fig.~6 keep the orientation and homological markings of traces.
If three traces marked by $[a],[b],[c]$ meet in a triple vertex 
 then $b=a+c\pmod{|m|}$, where $[b]$ is the marking of the middle trace,
 see a more general case in \cite[Lemma~6.3]{FK}.
The trace graph always remains symmetric under $\tau\mapsto\tau+\pi$.
Hence each move of Fig.~6 describes how to replace a small disk $\De$
 and the symmetric image of $\De$ under $\tau\mapsto\tau+\pi$ by 
 another small disk $\De'$ and the symmetric image of $\De'$, respectively.
\end{proof}

\begin{figure}[!h]
\includegraphics[scale=1.0]{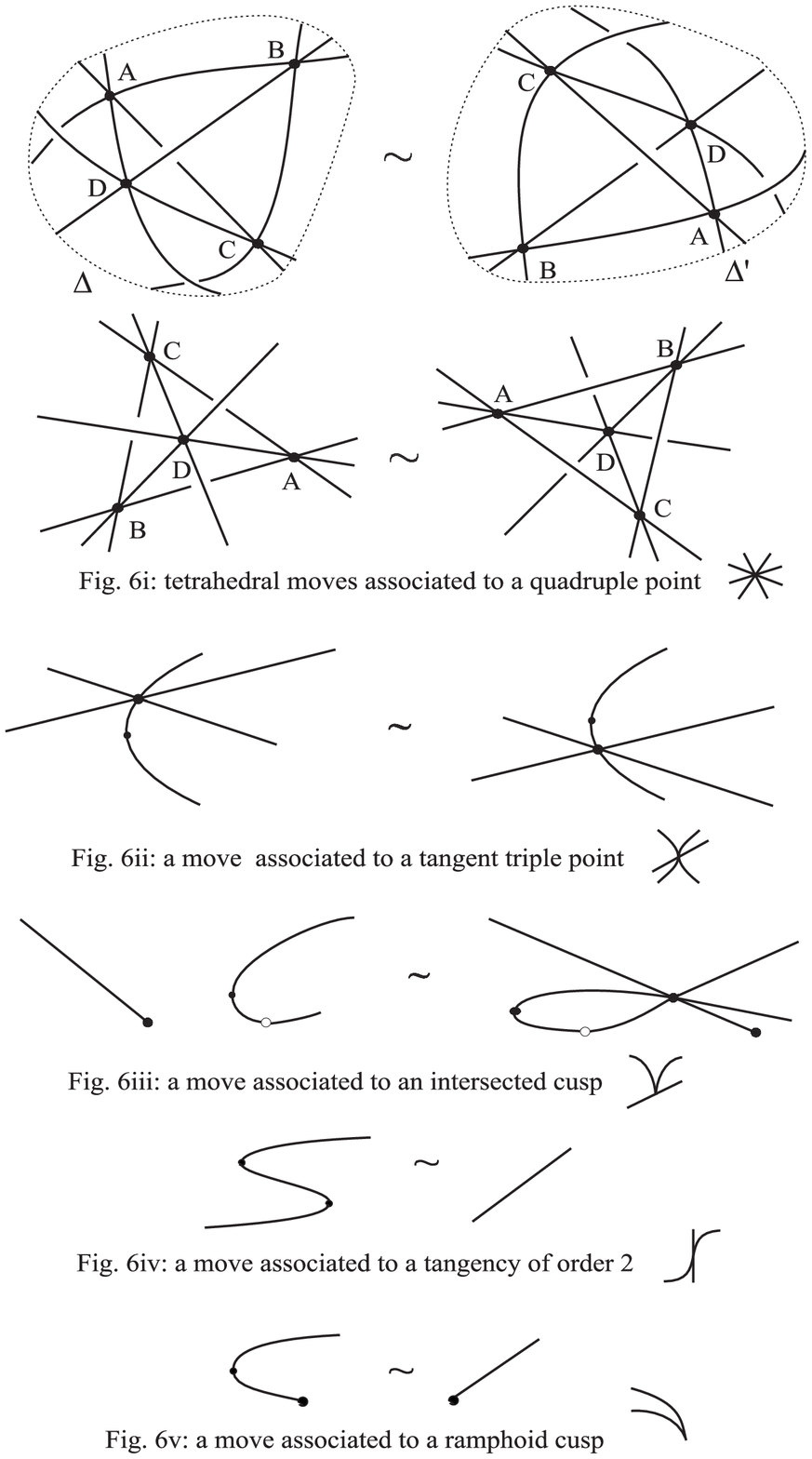}
\end{figure}

\begin{figure}[!h]
\includegraphics[scale=1.0]{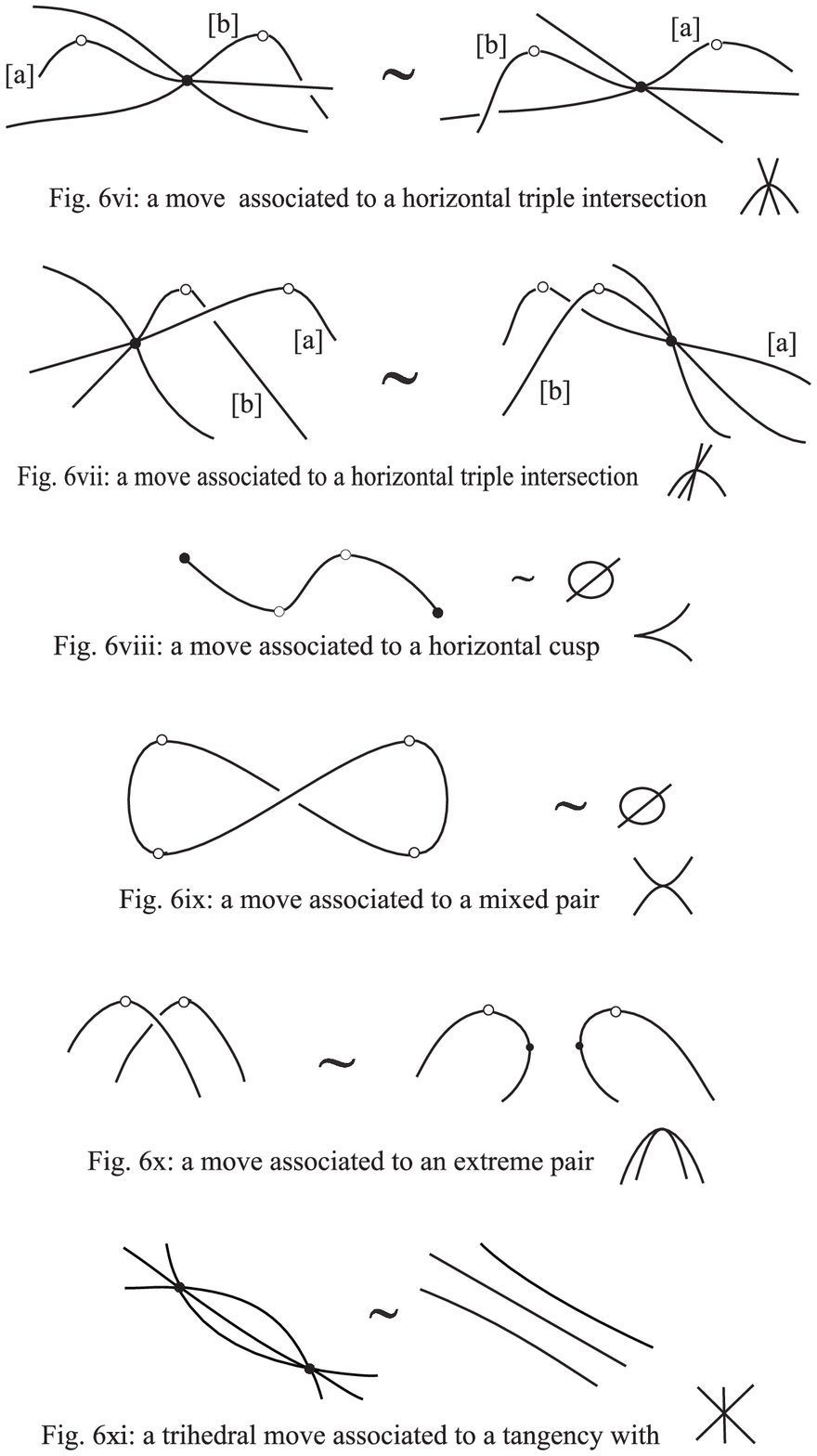}
\caption{Moves on trace graphs}
\end{figure}


\section{Proofs of main results}

\begin{lemma}
The writhes of Definition~1.3 are invariant under
 regular isotopy of trace graphs in the sense of Definition~2.5.
\end{lemma}
\begin{proof}
The Reidemeister moves of types II, III and IV in Fig.~4
 do not change the sum of signs in the writhes.
The Reidemeister move of type V either adds or deletes a crossing of $\tg(K)$,
 but a trace arc coming to a hanging vertex always has homology marking $[0]$
 modulo $|\lk(K,L)|$ and is excluded in Definition~1.3.
\end{proof}

\begin{lemma}
The moves in Fig.~6 keep the writhes of Definition~1.3 except
\smallskip

\noindent
$\bu$ 
the move 6i changes $W^u_{a,b}$ ($a\neq b$) by $\pm 2$
 for at most 6 unordered pairs $\{a,b\}$;
\smallskip

\noindent
$\bu$ 
the move 6i changes $W^o_{a,b}$ ($a\neq b$) by $\pm 1$
 for at most 12 ordered pairs $(a,b)$;
\smallskip

\noindent
$\bu$ 
the move 6i changes $W^c_{a,a}$  
 either (1) by $\pm 6$ for at most 2 values of $a$, or \\
\hspace*{2mm} 
(2) by $\pm 4$ for at most 2 values of $a$ 
 and by $\pm 2$ for at most 2 values of $a$, or \\
\hspace*{2mm} 
(3) by $\pm 2$ for at most 6 values of $a$;
\smallskip

\noindent
$\bu$ 
the moves 6vi, 6vii change $W^o_{a,b}$ ($a\neq b$) by $\pm 1$
 for at most 4 ordered pairs $(a,b)$;
\smallskip

\noindent
$\bu$ 
the moves 6ix and 6x change $W^c_{a,a}$ by $\pm 1$ for at most two values of $a$.
\end{lemma}
\begin{proof}
The move 6i switches exactly 3 couples of symmetric crossings. 
For instance, the arc $DB$ overcrosses $AC$ in the left picture of Fig.~6i,
 but $DB$ undercrosses $AC$ in the right picture.
Hence, for at most 6 unordered pairs $\{a,b\}$ with $a\neq b$,
 the unordered writhe $W^u_{a,b}$ changes by $\pm 2$.
Similarly, for at most 12 ordered pairs $(a,b)$, 
 the ordered writhe $W^o_{a,b}$ changes by $\pm 1$
 since exactly one crossing of a trace $[a]$ over a trace $[b]$ 
 either appears or disappears under the move 6i.
\smallskip

If all 3 crossings in the disk $\De$ in Fig.~6i
 are formed by traces with same homological marking $[a]$ then
 the coordinated writhes $W_{a,a}^c$ and $W_{|m|-a,|m|-a}^c$ change
 by $\pm 6$ as required in the case (1).
If two of the above crossings are formed by a trace $[a]$
 and the remaining one by a different trace $[b]$
 then we arrive at the case (2).
The case (3) arises when each of the 3 crossings in $\De$ is formed
 by a different trace.
\smallskip

In the moves 6vi and 6vii the overcrossing arc becomes undercrossing
 and vice versa, but the sign of the crossing is invariant,
 i.e. $W_{a,b}^u$ does not change.
Each of the moves 6vi and 6vii deletes exactly one crossing, 
 where a trace $[a]$ overcrosses a trace $[b]$,
 and adds another crossing, where the trace $[a]$ undercrosses the trace $[b]$.
Under the symmetry $\tau\mapsto\tau+\pi$, we get similar conclusions
 for the traces marked by $[|m|-a]$ and $[|m|-b]$.
So the ordered writhe $W_{a,b}^o$ changes by $\pm 1$
 for the 4 ordered pairs $(a,b)$, $(b,a)$ and $(|m|-a,|m|-b)$, $(|m|-b,|m|-a)$.
\smallskip
 
The move 6ix adds or deletes a crossing of 
 a trace circle $[a]$ with itself.
Hence only the writhes $W_{a,a}^c$ and $W^c_{|m|-a,|m|-a}$ change by $\pm 1$.
The move 6x adds or deletes a crossing between two arcs
 belonging to traces with same homological marking $[a]$.
Indeed, a pair of crossings corresponding to these arcs
 looks like in a horizontal version $\maxtang$ of Reidemeister move II, see Fig.~4. 
By Definition~2.2 the markings of these crossings are equal.
So the conclusion is the same as for the move 6ix.
\end{proof}
\medskip

\noindent
{\bf Proof of Theorem~1.4.}
To prove the first lower bound it suffices to show that
 the right hand side increases by 1 only if a generic isotopy $\{K_t\}$ 
 passes through a knot with a fiber quadrisecant.
By Lemma~3.2 the unordered writhe changes under the move 6i
 associated to a fiber quadrisecant, see the correspondence between
 singularities and moves in \cite[section~4]{FK}.
Six unordered pairs $\{a,b\}$ provide
 the maximal increase 1 as required.
Lemma~3.2 also proves the third lower bound since 
 only the moves 6i, 6ii, 6iv, 6xi are relevant for braids.
\smallskip

For the second lower bound, we are interested in
 crossings, where arcs have same marking.
By Lemma~3.2 the coordinated writhe changes only under the move
 6i and two moves 6ix, 6x associated to
 a fiber extreme secant in knot isotopies.
Under the move 6i the right hand side increases at most by 1
 while under the moves 6ix and 6x the maximal increase is 1/6 
 after multiplying by 1/12.
\qed
\medskip

\begin{figure}[!h]
\includegraphics[scale=1.0]{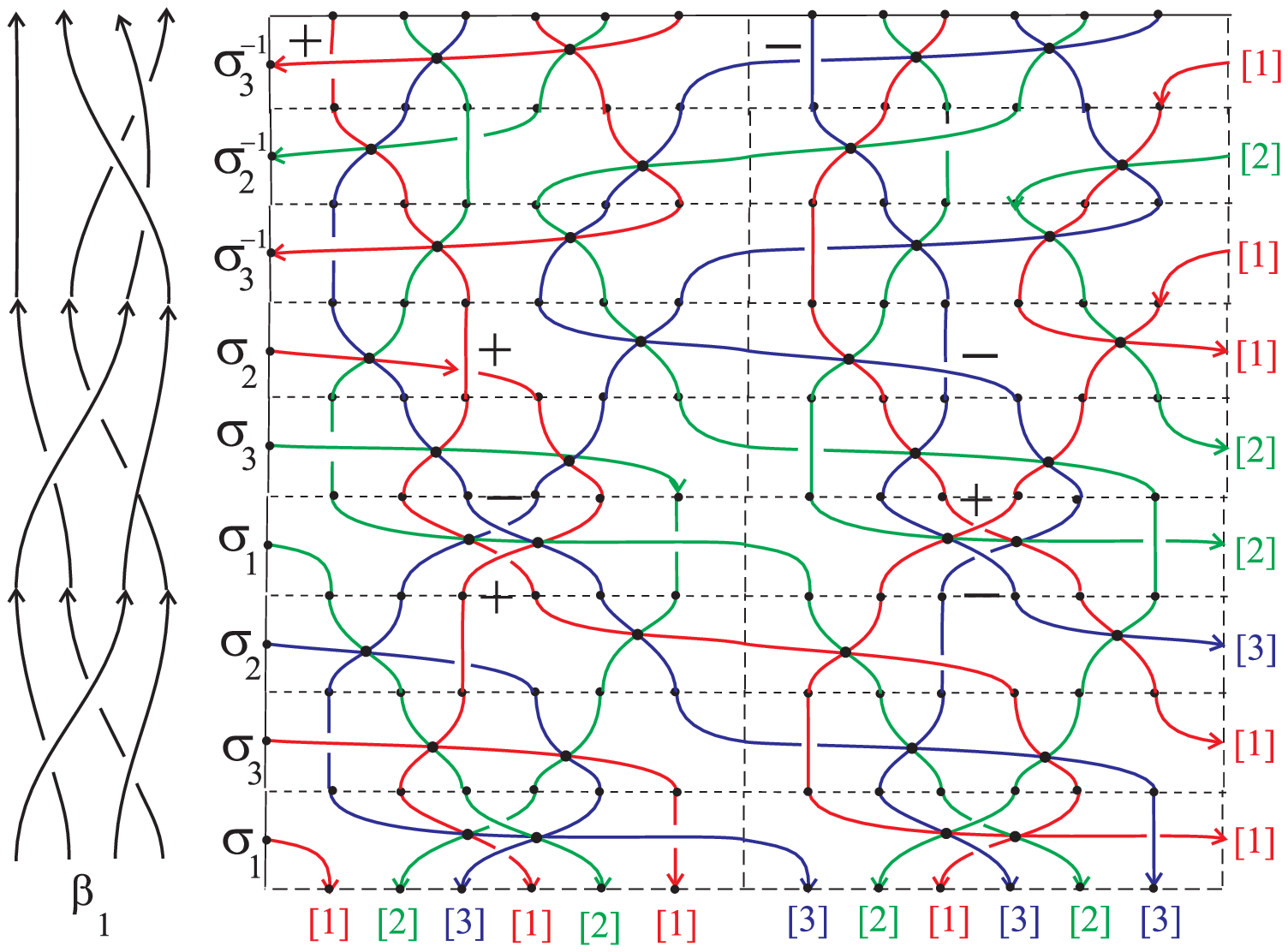}
\caption{The trace graph of the closure of
 $\be_1=(\si_1\si_3\si_2)^2\si_3^{-1}\si_2^{-1}\si_3^{-1}$.}
\end{figure}

\begin{example}
Consider the isotopic closures of the braids
 $\be_0=\si_3\si_2\si_1$ and $\be_1=(\si_1\si_3\si_2)^2\si_3^{-1}\si_2^{-1}\si_3^{-1}$.
The trace graphs of $\hat\be_0$ and $\hat\be_1$ are in Fig.~3,7, respectively.
They were constructed by attaching elementary blocks described in the proof
 of Lemma~2.4.
So we assume that the closed braids are given by embeddings
 into a neighbourhood of the torus $S_{\tau}^1\times S_{\ph}^1$
 located vertically in $\R^3-L$.
\smallskip

Both graphs split into 3 closed traces (circles with self-intersections)
 marked by $[1],[2],[3]$.
The trace graph in Fig.~3 has no crossings, i.e. 
 the writhes of Definition~1.3 vanish.
For the trace graph of $\hat\be_1$, the non-zero writhes 
 are $W_{1,1}^c=4$, $W^c_{3,3}=-4$. 
The 4 signs $+$ and 4 signs $-$ are shown in Fig.~7.
The second lower bound of Theorem~1.4 implies that any isotopy connecting 
 the closed braids $\hat\be_0,\hat\be_1$ involves at least one fiber quadrisecant.
The conclusion is the same for the closures of $\be_0\ga,\be_1\ga$,
 where $\ga$ is any pure 4-braid.
\smallskip

Consider the 4-braid $\be=(\si_1\si_3\si_2)^2(\si_1^{-1}\si_3^{-1}\si_2^{-1})^2$
 in Fig.~8 and the sequence of the braids $\be_n=\be^{n-1}\be_1$, $n\geq 1$,
 whose closures are isotopic to $\hat\be_0$, see Fig.~3.
Fig.~8 contains the part of $\tg(\hat\be_n)$ corresponding to 
 a single factor $\be$ in $\be_n$.
So $\tg(\hat\be_n)$ is obtained from $\tg(\hat\be_1)$ by inserting
 $n-1$ copies of Fig.~8 at the bottom of Fig.~7.
The part in Fig.~8 has the writhes $W^c_{1,1}=3$ and $W^c_{3,3}=-3$.
Hence $\tg(\hat\be_n)$ has $W^c_{1,1}=3n+1$ and $W^c_{3,3}=-3n-1$.
By Theorem~1.4 any isotopy connecting the closures of $\be_0$
 and $\be_n$ involves at least $\dfrac{3n+1}{12}$ fiber quadrisecants.
So the second lower bound of Theorem~1.4 can be arbitrarily large.
\end{example}

\begin{figure}[!h]
\includegraphics[scale=1.0]{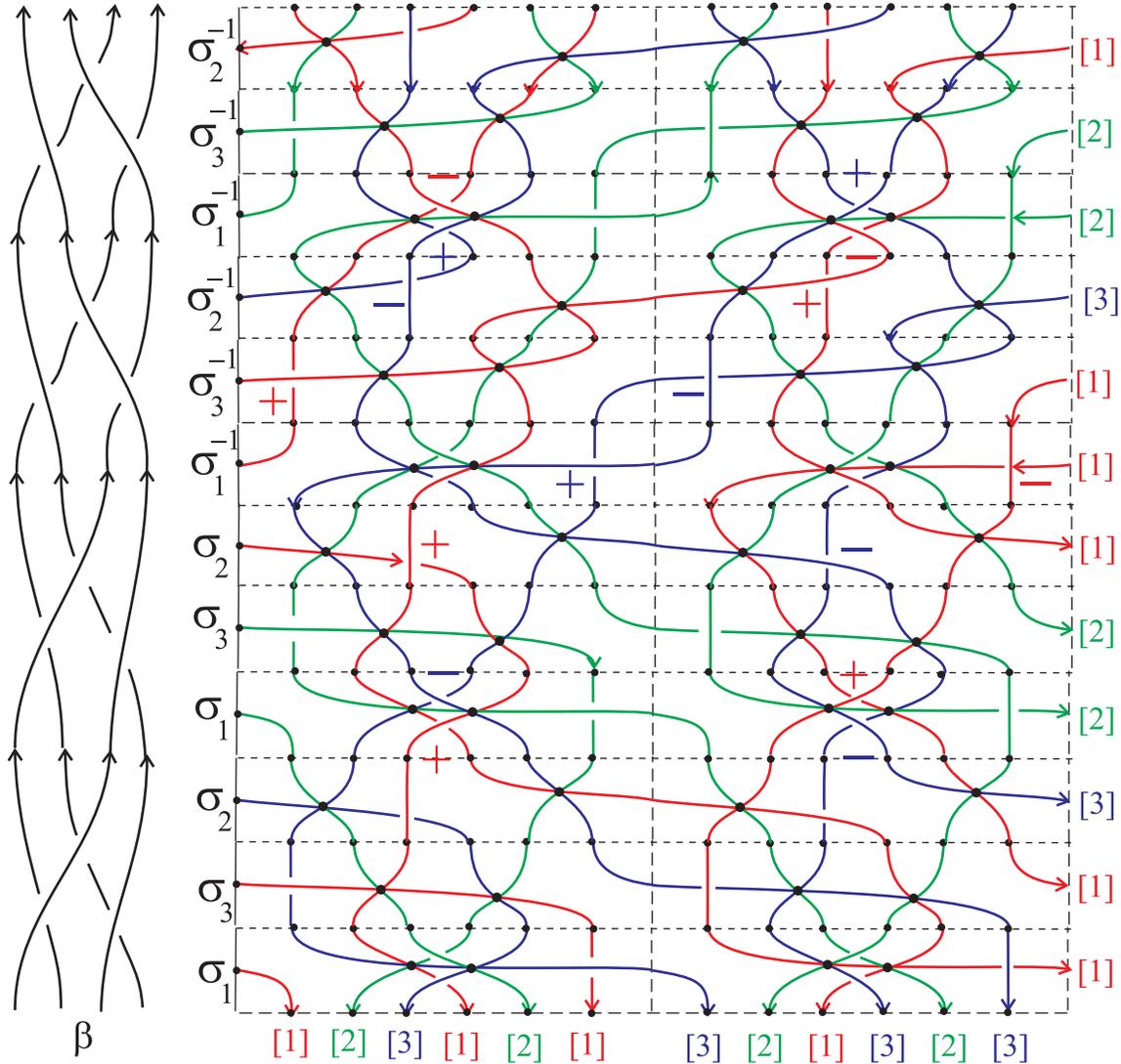}
\caption{The part of the trace graph for the factor 
 $\be=(\si_1\si_3\si_2)^2(\si_1^{-1}\si_3^{-1}\si_2^{-1})^2$.}
\end{figure}


\end{document}